\documentclass[letterpaper, 10 pt, conference]{ieeeconf}  
\usepackage{hhline}
\usepackage[letterpaper,
            left=.75in,
            right=.75in,
            top=1in,
            bottom=.75in]{geometry}
\usepackage{tikz}
\usetikzlibrary{positioning}
\usetikzlibrary{arrows.meta,positioning,fit}
\usepackage{cite}
\usepackage{amsmath,amssymb,amsfonts}
\usepackage{bm}
\usepackage{algpseudocode}
\usepackage{algorithm}
\usepackage{graphicx}
\usepackage{textcomp}
\usepackage[abs]{overpic}
\usepackage{mathrsfs}
\usepackage{xcolor}
\usepackage{hyperref}
\definecolor{lightblue}{rgb}{0.68, 0.85, 0.9}

\newcommand{\bbE}{\mathbb{E}}

\newcommand{\bbH}{\mathbb{H}}

\newcommand{\bbL}{\mathbb{L}}

\newcommand{\bbR}{\mathbb{R}}

\newcommand{\bbT}{\mathbb{T}}

\newcommand{\X}{\mathcal{X}}
\newcommand{\Z}{\mathcal{Z}}

\newcommand{\D}[3]{D_{#1}\!\left(#2,#3\right)}

\newcommand{\bx}{\bm{x}}

\newcommand{\rate}{e^{a_t + b_t}}
\newcommand{\Ssum}{\sum_{i=1}^N}

\newtheorem{theorem}{Theorem}[section]

\newtheorem{lemma}[theorem]{Lemma}

\newtheorem{definition}[theorem]{Definition}

\IEEEoverridecommandlockouts                              

\title{\LARGE \bf
Unifying Variational View of Accelerated Primal-Dual Methods
}

\author{Cheng Chang, \IEEEmembership{Members, IEEE} and Mehran Mesbahi \IEEEmembership{Fellow, IEEE}
\thanks{The research of the authors has been supported by the Air Force Office of Scientific Research grant FA9550-26-1-B206.}
\thanks{The authors are with the William E. Boeing Department of Aeronautics and Astronautics, University of Washington; emails: \{chang53+mesbahi\}@uw.edu.}}

\begin{document}

\maketitle
\thispagestyle{empty}
\pagestyle{empty}

\begin{abstract}

We develop a unifying variational framework for accelerated primal–dual flows in affinely constrained convex optimization. In particular, it is shown that applying the \emph{Euler–Lagrange–Rayleigh} (ELR) principle to coupled augmented Bregman Lagrangians systematically generates accelerated dynamics over general Bregman geometries and recovers a family of existing primal–dual mirror and Alternating Direction Method of Multipliers (ADMM) flows as special cases. We then establish $\mathcal{O}(e^{-b_t})$ convergence guarantees for these flows under mild assumptions. Lastly, we show that the proposed framework for algorithmic development extends from finite-dimensional optimization to constrained optimization over probability distributions.
\end{abstract}

\section{INTRODUCTION}

Continuous-time accelerated flows for optimization have received considerable attention following the elegant realization that the continuous limit of Nesterov's accelerated gradient descent corresponds to second-order ODEs~\cite{continuous_nestrov}. This correspondence has inspired many subsequent works on the differential equations-based analysis governing acceleration, including its variational interpretation proposed in \cite{BD_flow1} via the \emph{Euler-Lagrange} (EL) equation on the Bregman Divergence (BD) Lagrangian. Subsequently, \cite{dual_flows_1} constructed a dual BD Lagrangian for the affinely-constrained problem, although full exponential ($\mathcal{O}(e^{-b_t}),\;b_t:\bbT \to \bbR^+,\;b_t \to \infty$) rates still rely on a strong-convexity assumption. Meanwhile, several accelerated flows for affinely-constrained convex problems have been proposed independently, including accelerated second-order ODEs \cite{primal_dual_2,decentralized_flow2}, accelerated mirror-space flows \cite{decentralized_flows,primal_dual_1} under the standard scaling $b_t=\log K+2\log t$ ($K>0$), mixed-type flows with first-order dual ascent \cite{primal_dual_new}, and accelerated ADMM dynamics~\cite{admm_flows}. Other perspectives on acceleration have also been proposed; in \cite{feedback_control_optimization}, feedback control generates the multiplier dynamics, while \cite{amir_BD_flow} formulates the unconstrained variational problem as an optimal-control problem and obtains Hamiltonian flows through the Legendre transform. In the meantime, other types of acceleration methods have been examined using closely related constructions, hinting
at a common underlying unifying mechanism for their analysis.

In this paper, we identify such a unified perspective through a variational mechanism via the ELR equation \cite{Euler_Lagrange_Rayleigh} applied to coupled BD Lagrangians constructed from the augmented Lagrangian. In this formulation, the interaction between accelerated primal descent and multiplier ascent follows from the EL stationary-action construction extended by Rayleigh dissipation, while a common Lyapunov structure exposes the key convergence properties. Many existing methods \cite{admm_flows,decentralized_flows,primal_dual_2,decentralized_flow2} then arise as special cases of proposed variational framework, providing a systematic route to construct accelerated flows and establish their convergence rates. Furthermore, unlike formulations restricted to Euclidean dynamics, we show that the BD construction permits broader local geometries generated by strongly convex functions, allowing the primal and multiplier dynamics to evolve through mirror coordinates, including the negative-entropy geometry with \emph{Kullback--Leibler} (KL) divergence. We also extend the unconstrained distributional accelerated flows in \cite{amir_BD_flow} to affinely-constrained optimization over probability distributions within the same ELR framework. The main contributions of this paper are as follows: (i) providing a unified ELR interpretation of accelerated primal-dual flows; (ii) introducing a geometry-flexible construction admitting general BD generators beyond the Euclidean case; and (iii) providing an extension to constrained distributional optimization with accelerated dual-gap, feasibility, and objective convergence guarantees under the stated assumptions. Thus, this paper does not merely propose another acceleration dynamics, but unifies current primal-dual acceleration dynamics under the proposed ELR perspective. 
Lastly, to the best of our knowledge, this is the first variational-based accelerated framework for the constrained distributional setting.
In this direction, Table~\ref{table: comparison} summarizes the representative methods. 

\begin{table}[H]
\caption{Summary of representative existing and current works}
\label{table: comparison}
\begin{center}
\begin{tabular}{ccccc}
  & Feasibility & Objective& Dual gap& Variational view
\\
\hhline{=====}
This work & $\mathcal{O}(e^{-b_t})$ & $\mathcal{O}(e^{-b_t})$ & $\mathcal{O}(e^{-b_t})$ & Coupled ELR \\

\cite{BD_flow1} &None &$\mathcal{O}(e^{-b_t})$ &None& Primal EL\\

\cite{dual_flows_1} &$\mathcal{O}(e^{-b_t/2})$&STR,$\mathcal{O}(e^{-b_t})$&$\mathcal{O}(e^{-b_t})$& Dual EL\\

\cite{feedback_control_optimization} &$\mathcal{O}(e^{-b_t/2})$&STR,$\mathcal{O}(e^{-b_t})$&$\mathcal{O}(e^{-b_t})$& None\\

\cite{admm_flows}&$\mathcal{O}(e^{-b_t/2})$&$\mathcal{O}(e^{-b_t})$&$\mathcal{O}(e^{-b_t})$&None\\

\cite{decentralized_flows,primal_dual_2}&$\mathcal{O}(\frac{1}{t^2})$&$\mathcal{O}(\frac{1}{t^2})$& $\mathcal{O}(\frac{1}{t^2})$ & None
\\
\cite{primal_dual_new} & $\mathcal{O}(e^{-b_t})$ & $\mathcal{O}(e^{-b_t})$ & $\mathcal{O}(e^{-b_t})$ & None \\
\cite{primal_dual_1,decentralized_flow2}&$\mathcal{O}(\frac{1}{t})$&$\mathcal{O}(\frac{1}{t})$&$\mathcal{O}(\frac{1}{t^2})$&None\\
\hhline{=====}
\multicolumn{5}{p{242pt}}{STR: Strongly convex}
\end{tabular}
\end{center}
\end{table}

\section{Problem Formulation}
Consider the following jointly constrained multi-variable optimization problem on a finite-dimensional vector space with \emph{convex} and \emph{continuously differentiable} objectives $f_i:\X_i\to \bbR$: 
\begin{equation} \label{problem 1}
      \min_{x_i\in\X_i, \forall i\in \mathcal{N}}\ \Ssum f_i(x_i)
  \quad\text{s.t.}\quad
  \Ssum A_i x_i=C,
\end{equation}
where $\X_i\subseteq \bbR^n$ is convex and nonempty, $C\in \Z\subseteq\bbR^m$, $A_i\in \bbR^{m\times n}$, and $\mathcal{N} = \{1, \ldots,N\}$ is the index set. Assume that \eqref{problem 1} is feasible, admitting an optimal pair $(\bx^*,\lambda^*)$ satisfying the Karush–Kuhn–Tucker (KKT) conditions.
In this work, we initially let $\X_i = \bbR^n$ for simplicity. Note that certain convex set constraints can be posed by carefully designing the mirror mapping in Fig.\ref{fig:coordinate-map} \cite{decentralized_flows}. A numerical example considering the positive simplex set is given in \S\ref{subsec: example 1}.

Let $\bx = (x_1,\ldots,x_N)\in \X:= \prod_{i=1}^N \X_i$ and their corresponding feasible sets with objective function $F(\bx) = \Ssum f_i(x_i)$; we can then define linear operator $\mathcal{A}\bx:= \Ssum A_ix_i$.The feasibility can then be conveniently expressed in terms of the residue $r(\bx) = \mathcal{A}\bx - C$. 
This formulation represents distinct classes of optimization algorithms, including the standard ADMM when $N=2$ \cite{admm_flows,admm} and certain distributed optimization \cite{admm,admm_flows,decentralized_flows}. In the next section, we construct the augmented BD Lagrangians with the EL condition to reveal the underlying accelerated flows.

\section{ACCELERATED PRIMAL-DUAL FLOWS} \label{sec: primal dual flows}
Inspired by the BD Lagrangian for the unconstrained problem first mentioned in \cite{BD_flow1}, we will define the BD Lagrangian with the augmented Lagrangian for the constrained optimization defined in \eqref{problem 1}. Let $\mu>0$, the augmented Lagrangian for \eqref{problem 1} is 
\begin{equation} \label{eq: augmented lagrangian}
    \mathcal{L}_{\mu}(\bx,\lambda) = F(\bx) + \langle\lambda, r(\bx)\rangle + \frac{\mu}{2}\|r(\bx)\|^2,
\end{equation}
where $\lambda \in \bbR^m$ is the multiplier associated with the constraint. To this end, we formally define the BD and augmented BD Lagrangian as follows:

\begin{definition}[Bregman divergence]\label{def:BD}
Bregman divergence (BD) is a divergence equipped with a \emph{strictly convex} and \emph{differentiable} function $\phi_i: \X_i \to \bbR$ \cite{BD_lqr} and is written as:
\[
    D_{\phi_i}(y,x) =  \phi_i(y) - \phi_i(x) - \langle\nabla\phi_i(x), y - x\rangle,\quad \forall x,y \in \X_i.
\]
BD has several useful properties and examples that make it a more general measurement of discrepancy, though it lacks the common axioms of a metric, such as symmetry and the triangle inequality. To see more examples of the BD, refer to \cite{amir_BD_flow} for the KL divergence, and \cite{BD_lqr} for general properties.

\end{definition}

\begin{definition}[Augmented BD Lagrangians for \eqref{problem 1}]\label{def:BD lag}
Assume the \emph{strongly convex} and twice differentiable distance generating functions $\phi_i: \X_i \to \bbR$ and $\Phi_{\lambda}: \bbR^m \to \bbR$. 
Choose time-varying weighing functions $a_t,\;b_t,\;g_t: \bbT\to \bbR$ 
that satisfy the ideal scaling \cite{BD_flow1} as follows:
\begin{equation} \label{eq: ideal scaling}
    \dot{b}_t \leq e^{a_t},\quad \dot{g}_t = e^{a_t}.
\end{equation}
Let the lookahead variable be $z_{\bx} = \bx + e^{-a_t}\dot\bx = (z_{x_1},\ldots,z_{x_N})$ and $z_{\lambda} = \lambda + e^{-a_t}\dot\lambda$. One can think of the lookahead terms as producing second-order information in time.
Now define the Augmented BD Lagrangians for the primal variable and the multiplier with other variables fixed \emph{under variation}, e.g., fix $z_{\lambda}$ in \eqref{eq: primal augmented lagrangians}.
\begin{subequations} \label{eq: augmented lagrangians}
\begin{align} \label{eq: primal augmented lagrangians}
        \bbL_{\bx}(\bx,\dot\bx,t) &= e^{a_t+g_t}[\underbrace{D_{\Phi_x}(z_{\bx},\bx)}_{\text{Kinetic term}} - \underbrace{e^{b_t}\mathcal{L}_{\mu}(\bx,z_{\lambda})}_{\text{Potential term}}] \\
            \bbL_{\lambda}(\lambda,\dot\lambda,t) &= e^{a_t+g_t}[\underbrace{D_{\Phi_{\lambda}}(z_\lambda,\lambda)}_{\text{Kinetic term}} + \underbrace{e^{b_t}\mathcal{L}_{\mu}(z_{\bx},\lambda)}_{\text{Potential term}}],
\end{align}
\end{subequations} 
where $\Phi_x(\bx)= \Ssum\phi_i(x_i)$, $D_{\Phi_x}(p,q) = \Ssum D_{\phi_i}(p_i,q_i),\quad \forall p,q \in \X$, and $\nabla \Phi_{\bx} = (\nabla\phi_1(x_1),\ldots,\nabla\phi_N(x_N))$. 

\end{definition}

Unlike the unconstrained case where the objective acts as a potential, the augmented Lagrangian \eqref{eq: augmented lagrangians} plays the role of the potential, with the primal variable and multiplier, respectively, minimizing and maximizing this potential to close the \emph{duality gap}. 

\subsection{Undamped ODE flows from Euler-Lagrange condition} \label{subsec: undamped flows}

Define the action functional $\mathcal{S}(q|_{t_0}^{t_f}) = \int_{t_0}^{t_f}\mathbb{L}_q(q,\dot q,t)dt$ \cite{calculus_of_variatino_optimal_control}, where $q \in \bbR^d$ denotes either $x_i$ $(d = n)$ or $\lambda$ ($d = m$), and $q^*|_{t_0}^{t_f}$ is the stationary path.
For least action trajectories induced by Lagrangian mechanics, from the calculus of variations, the locally optimal or stationary action paths (\emph{extremals})  must satisfy the \emph{Euler-Lagrange} condition for all admissible variations $v|_{t_0}^{t_f},\;\text{s.t. }v(t_0) = v(t_f) = 0$ \cite{calculus_of_variatino_optimal_control}:
\begin{equation} \label{eq: euler lagrangian}
\begin{aligned}
    q^*|_{t_0}^{t_f} \Leftrightarrow D\mathcal{S}(q^*|_{t_0}^{t_f})[v|_{t_0}^{t_f}] = 0 \Leftrightarrow 
    \frac{d}{dt}\frac{\partial\bbL_q}{\partial \dot q} - \frac{\partial \bbL_q}{\partial q} = 0.
\end{aligned}
\end{equation}\cite{calculus_of_variatino_optimal_control}. By solving \eqref{eq: euler lagrangian}, one can obtain sets of second-order ODEs for both primal variables $x_i$ and the multiplier $\lambda$.
The lemma below provides the corresponding first- and second-order ODEs for $x_i$ and $\lambda$.

\begin{lemma}\label{lemma:Undamped ODES}
Assume the ideal scaling \eqref{eq: ideal scaling} holds, then 
given the Euler-Lagrangian condition \eqref{eq: euler lagrangian} and the Augmented BD Lagrangians given in \eqref{eq: augmented lagrangians}, the resulting second-order ODEs for the primal variable and multipliers are:
\begin{subequations} \label{eq: undamped second order flows}
\begin{align} 
    \ddot x_i + (e^{a_t}- \dot a_t)\dot x_i + e^{2a_t + b_t}[\nabla^2 \phi_i(z_i)]^{-1} \nabla_{x_i}\mathcal{L_{\mu}}(\bx, z_{\lambda}) = 0,\\
    \ddot \lambda + (e^{a_t}- \dot a_t)\dot \lambda - e^{2a_t + b_t}[\nabla^2 \Phi_{\lambda}(z_{\lambda})]^{-1} \nabla_{\lambda}\mathcal{L_{\mu}}(z_{\bx}, \lambda) = 0.
\end{align}
\end{subequations}
From these second-order ODEs, one can construct first-order ODEs for primal variables as 
\begin{subequations} \label{eq: primal undamped flows}
\begin{align}  \label{eq: primal undamped flows: position}
&\underbrace{\dot x_i = e^{a_t}[\nabla \phi^*_i(u_i) - x_i]}_{\text{Primal position flow}},\\ 
\label{eq: primal undamped flows: mirror}
&\underbrace{\dot u_i = -e^{a_t+b_t}  \nabla_{x_i}\mathcal{L_{\mu}}(\bx, z_{\lambda})}_{\text{Dual space  \emph{descent} mirror flow}} ,
\end{align}
\end{subequations}, where $u_i = \nabla \phi_i (z_i)$, $ z_i= x_i + e^{-a_t}\dot x_i$, and $*$ denotes the conjugate. We also denote $u_{\bx} = (u_1,\ldots,u_N)$ for the compact reference. 
Similarly, the first-order ODEs for the multiplier can be formulated as:
\begin{subequations}  \label{eq: multiplier undamped flow}
\begin{align}
\label{eq: multiplier undamped flow: position}
\underbrace{\dot \lambda = e^{a_t}[\nabla \Phi^*_{\lambda}(u_{\lambda}) - \lambda]}_{\text{Multiplier position flow}},\\
\label{eq: multiplier undamped flow: mirror}
\underbrace{\dot u_{\lambda} = e^{a_t+b_t}  \nabla_{\lambda}\mathcal{L_{\mu}}(z_{\bx}, \lambda) }_{\text{Dual space \emph{ascent} mirror flow}},
\end{align}
\end{subequations}, where $u_{\lambda} = \nabla \Phi_{\lambda} (z_\lambda)$ and $ z_{\lambda}= \lambda + e^{-a_t}\dot \lambda$.

\begin{proof}
    See Appendix \ref{proof: lemma Undamped} for details.
\end{proof}

\end{lemma}

\begin{figure}[H] 
\centering
\begin{tikzpicture}[
  node distance=0.5cm and 0.1cm, 
  box/.style={draw, rounded corners, minimum width=1.4cm, minimum height=1.0cm, align=center, thick},
  arr/.style={-{Latex[length=1mm]}, thick}
]

\node[box, fill=lightblue, draw=black] (mirror) { Dual (covector) space \\$(u_{\bx}, u_{\lambda})$\eqref{eq: primal undamped flows: mirror}, \eqref{eq: multiplier undamped flow: mirror}};

\node[box, fill=lightblue, draw=black, above=of mirror] (primal) {Vector space\\$(\bx, \lambda)$ \eqref{eq: primal undamped flows: position}, \eqref{eq: multiplier undamped flow: position}};

\draw[arr, black] ([xshift=-12pt]primal.south) -- node[midway, left=4pt] {$\nabla \Phi$} ([xshift=-12pt]mirror.north);

\draw[arr, black] ([xshift=12pt]mirror.north) -- node[midway, right=2pt] {$\nabla \Phi^*$} ([xshift=12pt]primal.south);

\end{tikzpicture}
\caption{Dual space mirror flow and mirror maps $\nabla \Phi$ and $\nabla \Phi^*$.}
\label{fig:coordinate-map}
\end{figure}
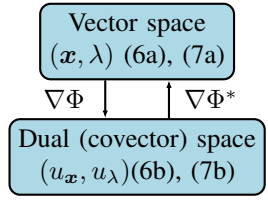

By separating the $N+1$ second-order ODEs in \eqref{eq: undamped second order flows} into $2(N+1)$ first-order ODEs in \eqref{eq: primal undamped flows} and \eqref{eq: multiplier undamped flow}, one can easily perform various numerical integrations to solve for the ODEs. Furthermore, as shown in Fig. \ref{fig:coordinate-map}, the position/ mirror flows can be related via the gradient of the chosen metric functions through $\nabla \Phi$. The descent and ascent are performed in the dual space before being mapped back to the corresponding vector spaces, resembling mirror descent methods \cite{mirror,accelerated_mirror} but performing descent and ascent simultaneously to close the dual gap.

\subsection{Pontryagin’s Maximum Principle and optimal control derivation}
This section provides an optimal control perspective formulation using the PMP through the construction of Hamiltonian flows \cite{calculus_of_variatino_optimal_control}.
From the Augmented BD Lagrangians defined earlier \eqref{eq: augmented lagrangians}, 
consider the following two optimal control problems with the velocity constraints and costs of Lagrangians:
\begin{subequations} \label{eq: optimal control problems}
\begin{align}
    \min_{v_{\bx}} &\int_0^{\infty} \bbL_{\bx} (\bx,\dot \bx,t) dt, \quad\text{s.t.} \quad \dot \bx = v_{\bx}, \\
    \min_{v_{\lambda}}& \int_0^{\infty} \bbL_{\lambda} (\lambda,\dot \lambda,t) dt, \quad\text{s.t.} \quad \dot \lambda = v_{\lambda},
\end{align}
\end{subequations}
where $v_{\bx} = \dot x$ and $v_{\lambda} = \dot \lambda$ are the velocity of the positions. These problems produce the corresponding Hamiltonians for $x$ and $\lambda$ with costates $Y_{\bx}$ and $Y_{\lambda}$:
\begin{subequations}
\begin{align}
    \bbH_{\bx}(\bx,v_{\bx},Y_{\bx},t) &= \langle Y_{\bx},v_{\bx}  \rangle - \bbL_{\bx}(\bx,v_{\bx}, t), \\
    \bbH_{\lambda}(\lambda,v_{\lambda},Y_{\lambda},t) &= \langle Y_{\lambda},v_{\lambda}  \rangle - \bbL_{\lambda}(\lambda,v_{\lambda}, t). 
\end{align}
\end{subequations} 

By Pontryagin’s Maximum Principle (PMP), the optimal velocity associated with these Hamiltonians is:
\begin{equation}
\begin{aligned}
    v_{\bx}^* &= \dot \bx = \frac{\partial \bbH_{\bx}}{\partial Y_{\bx}}\\
    & =\text{argmax}_{v_{\bx}} \bbH_{\bx}(\bx,Y_{\bx},t) = e^{a_t}[\nabla \Phi_{\bx}^* (u_{\bx})-\bx], \\
    v_{\lambda}^* &= \dot \lambda = \frac{\partial \bbH_{\lambda}}{\partial Y_{\lambda}}\\
    & = \text{argmax}_{v_{\lambda}} \bbH_{\lambda}(\lambda,Y_{\lambda},t) = e^{a_t}[\nabla \Phi_{\lambda}^* (u_{\lambda})-\lambda],
\end{aligned}
\end{equation} where we used the zero gradient conditions at the maximum: $0 = Y_q - e^{g_t}[\nabla \Phi_q(z_q) - \nabla\Phi_q(q)]$, where $q$ represents either $\bx$ or $\lambda$. Thus, the positional flows resulting from the optimal control problems \eqref{eq: optimal control problems} are the same as the Lagrangian formulations \eqref{eq: primal undamped flows: position} and \eqref{eq: multiplier undamped flow: position}. The costate (momentum) flows $\dot Y_{\bx} = -\frac{\partial\bbH_{\bx}}{\partial \bx}$ and  $\dot Y_{\lambda} = -\frac{\partial\bbH_{\lambda}}{\partial \lambda}$ can also be shown to directly connected to the mirror flows in \eqref{eq: primal undamped flows: mirror} and \eqref{eq: multiplier undamped flow: mirror} by applying $Y_q = \frac{\partial \bbH_q}{\partial \dot{q}}$ \cite{calculus_of_variatino_optimal_control}.

\subsection{Friction, Acceleration, and Convergence}
The EL in \eqref{eq: euler lagrangian} can be extended to include the velocity-dependent generalized Rayleigh \emph{dissipation} function $\mathcal{R}_{\tiny \mbox{fri}}(q,\dot q,t) = \frac{\zeta}{2}e^{g_t}\|\dot q\|_{\nabla^2 h_q(q)}^2,\; \zeta \geq 0$ with strongly convex and twice differentiable $h_q:\bbR^d \to \bbR$ ($q = x_i \text{ or } \lambda$) to account for the non-conservative friction force as:
\begin{equation}
\label{eq: euler lagrangian rayleigh}
\begin{aligned}
    D\mathcal{S}(q|_{t_0}^{t_f})[v|_{t_0}^{t_f}] + \int_{t_0}^{t_f} \langle \mathcal{F}_{\tiny \mbox{fri}} , v \rangle dt = 0 \Leftrightarrow \\
     \frac{d}{dt}\frac{\partial\bbL_q}{\partial \dot q} - \frac{\partial \bbL_q}{\partial q}  - \mathcal{F}_{\tiny \mbox{fri}}= 0, 
\end{aligned}
\end{equation}
with the velocity-dependent friction force $\mathcal{F}_{\tiny \mbox{fri}} = -\frac{\partial \mathcal{R}_{\tiny \mbox{fri}}}{\partial \dot q} = - \zeta e^{g_t} \nabla^2 h(q)\dot q$. In the mechanical system analogy, the Rayleigh dissipation term acts as damping to remove excess energy from the system and improve convergence stability; the resulting accelerated energy-decay bound can be shown rigorously via Lyapunov analysis. 
Through a similar derivation in \ref{proof: lemma Undamped}, the resulting second-order ODEs will thus change slightly and take the following forms:
\begin{equation} \label{eq:primal damped second order flows}
\begin{aligned} 
    \ddot x_i + (e^{a_t}- \dot a_t)\dot x_i + e^{2a_t + b_t}[\nabla^2 \phi_i(z_i)]^{-1} \nabla_{x_i}\mathcal{L_{\mu}}(\bx, z_{\lambda}) \\+ \zeta e^{a_t}[\nabla^2 \phi_i(z_i)]^{-1} \nabla^2h_i(x_i)\dot x_i= 0,
\end{aligned}
\end{equation}
and
\begin{equation} \label{eq:multiplier damped second order flows}
\begin{aligned} 
    \ddot \lambda + (e^{a_t}- \dot a_t)\dot \lambda - e^{2a_t + b_t}[\nabla^2 \Phi_{\lambda}(z_{\lambda})]^{-1} \nabla_{\lambda}\mathcal{L_{\mu}}(z_{\bx}, \lambda)\\ + \zeta e^{a_t}[\nabla^2 \Phi_{\lambda}(z_{\lambda})]^{-1} \nabla^2 h_{\lambda}(\lambda)  \dot \lambda= 0.
\end{aligned}
\end{equation}
Similarly, the mirror flows of $\bx$ and $\lambda$ now takes the form of 
\begin{subequations} \label{eq: damped mirror flows}
\begin{align}
    &\dot u_i = -e^{a_t+b_t}  \nabla_{x_i}\mathcal{L_{\mu}}(\bx, z_{\lambda}) - \zeta \nabla^2h_i(x_i) \dot x_i, \\
    &\dot u_{\lambda} = e^{a_t+b_t}  \nabla_{\lambda}\mathcal{L_{\mu}}(z_{\bx}, \lambda) - \zeta \nabla^2h_{\lambda}(\lambda) \dot \lambda,
\end{align}
\end{subequations}
for some $\zeta\geq 0$.

We aim to show that under the assumptions given before, the dual gap from the ODEs \eqref{eq: primal undamped flows},\eqref{eq: multiplier undamped flow} and damped mirror ODEs \eqref{eq: damped mirror flows} will lead to the \emph{exponential} rates of several properties, including the dual gap, constraint residue, and the objective function.

\begin{theorem}[ $\mathcal{O}(e^{-b_t})$ convergence rates]\label{theorem: convergence} Define the dual gap as: 
\begin{equation} \label{eq: dual gap}
    G_{\mu}(\bx, \lambda) = \mathcal{L}_{\mu}(\bx,\lambda^*) - \mathcal{L}_{\mu}(\bx^*,\lambda) = \mathcal{L}_{\mu}(\bx,\lambda^*) - F(\bx^*).
\end{equation} 

We claim that the following convergence results hold under the assumptions made earlier.
\begin{itemize}
  \item $G_{\mu}(\bx,\lambda) = \mathcal{O}(e^{-b_t})$. 
\item $\|r(\bx)\|^2 = \mathcal{O}(e^{-b_t})$. 
  \item $|F(\bx) - F(\bx^*)|^2 = \mathcal{O}(e^{-b_t})$.
\end{itemize}
Furthermore, if we assume that $\|u_{\lambda} + \zeta \nabla h_{\lambda}(\lambda) \| \leq K< \infty$, we have ``full'' decay rates $\mathcal{O}(e^{-b_t})$ for $\|r(\bx)\|$ and $|F(\bx) - F(\bx^*)|$. A sharper energy dissipation bound is obtained if $\zeta >0$, acting as extra damping. The boundedness holds when $\Phi_{\lambda}$ and $h_{\lambda}$ are strongly convex (Lemma \ref{lemma: boundedness}).

\begin{proof}
The core mechanism lies in the construction of Lyapunov functions composed of mixed measures of kinetic energy, dual gap, and distance to optimal as follows:
\begin{equation}
\begin{aligned} \label{eq:energy}
  \mathcal{E}=&
  \underbrace{\Ssum \D{\phi_i}{x_i^*}{z_i}}_{1. \text{Primal kinetic term}}
  +\underbrace{\D{\Phi_\lambda}{\lambda^*}{z_\lambda}}_{2. \text{Multiplier kinetic term}}
  +\underbrace{e^{b_t}G_\mu(\bx,\lambda)}_{3. \text{Potential term}}\\
  &+\underbrace{\zeta
  \left[
    \Ssum D_{h_i}(x^*_i,x_i)
    + D_{h_{\lambda}}(\lambda^*,\lambda)
  \right]}_{4. \text{Distance term}},
\end{aligned}
\end{equation}
where $x^*_i \in \X_i$ and $\lambda^*$ satisfy the KKT condition, e.g., zero-dual gap, and both primal and dual feasibility. The detailed proof is given in \ref{proof: convergence theorem}.
\end{proof}

\end{theorem}

\section{SPECIALIZATION TO UNCOVER EXISTING METHODS}
Using the Augmented BD Lagrangian formulation, we have identified existing classes of accelerated algorithms for constrained optimization as special instances of our proposed framework based on the fundamental ELR condition \eqref{eq: euler lagrangian rayleigh}. As such, the underlying ODE flows of many separate accelerated constrained optimization methods result from applying the variational principle to artificially constructed mechanical systems, with the augmented Lagrangian \eqref{eq: augmented lagrangian} as the system potential.

Next, we select representative specializations of this framework to demonstrate the generality of the proposed perspective. For the ease of comparison, we will preserve some of the original notation and expressions from the source papers.

\subsection{Accelerated primal-dual flows \cite{decentralized_flows}}
Consider the standard affinely constrained optimization problem in \cite{decentralized_flows}:
\begin{equation}\label{problem: primal dual}
\min_{x\in\X} f(x) \quad \text{s.t.} \quad Ax - B = 0,
\end{equation} where $\X \subseteq \bbR^n$, $A\in \bbR^{m \times n}$, and $B  \in \bbR^m$. The continuous accelerated flows on the primal-dual proposed by \cite{decentralized_flows} form the four ODEs:
\begin{equation}
\begin{aligned}  \label{eq: paper accelerated primal dual flow}
    \dot x &= \frac{\alpha}{t}(\nabla\Phi^*(u)-x),\quad \dot \lambda = \frac{\alpha}{t}(v-\lambda),\\
    \dot u &= -\frac{t}{\alpha}\underbrace{[\nabla f(x)+\mu A^\top(Ax - B) + A^\top z_{\lambda}]}_{\nabla_x\mathcal{L}(x,z_{\lambda})} - \zeta \dot x,\\
    \dot v &= \frac{t}{\alpha}\underbrace{(A\nabla\Phi^*(u)-B)}_{\nabla_{\lambda}\mathcal{L}(z_x,\lambda)} - \zeta \dot \lambda.
\end{aligned}
\end{equation}
Let $N=1$, $A_1 = A$, $C=B$, $\Phi_{\lambda} = \frac{1}{2}\|\cdot\|^2$, $h_{\lambda} = h_x = \frac{1}{2}\|\cdot\|^2$, and choose the following scalings: $a_t = \log\frac{\alpha}{t},$ $\quad b_t = 2\log\frac{t}{\alpha},$ and $ g_t = \alpha \log t$ for some $\alpha \geq 2$, then the ODE flows in \eqref{eq: primal undamped flows: position},\eqref{eq: multiplier undamped flow: position}, and \eqref{eq: damped mirror flows} can recover the accelerated primal-dual flows \eqref{eq: paper accelerated primal dual flow} proposed in \cite{decentralized_flows}.

\subsection{Accelerated ADMM\cite{admm_flows}}
Consider the standard ADMM problem posed in \cite{admm_flows}:
\begin{equation}\label{problem: admm}
\min_{x\in\X, y \in \mathcal{Y}} f(x) + g(y) \quad \text{s.t.} \quad Ax + By = C.
\end{equation} with $\X,\mathcal{Y} \subseteq \bbR^n$ and $A,B \in \bbR^{m \times n}$.
The accelerated ADMM flows proposed in \cite{admm_flows} for solving \eqref{problem: admm} take the form of 
\begin{subequations} \label{eq: admm flows}
\begin{align}
    &\ddot x + \gamma_t \dot x + \beta_t \nabla_x \mathcal{L}_{\mu}(x,y,\underbrace{\lambda +\alpha_t \dot \lambda}_{z_\lambda}) = 0, \\
    &\ddot y + \gamma_t \dot y + \beta_t \nabla_y \mathcal{L}_{\mu}(x,y,\underbrace{\lambda +\alpha_t \dot \lambda}_{z_\lambda}) = 0, \\
    &\ddot \lambda + \gamma_t \dot \lambda - \beta_t \nabla_{\lambda} \mathcal{L}_{\mu}(\underbrace{x + \alpha_t \dot x}_{z_x},\underbrace{y + \alpha_t \dot y}_{z_y},\lambda) = 0,
\end{align}
\end{subequations} with time-varying functions $\alpha_t,\beta_t,\gamma_t:\bbT\to \bbR^+$.

Let $N =2$, $x_1 = x$, $x_2 = y$, $\phi_1 = \phi_2 = \Phi_{\lambda} = h_1 = h_2=h_{\lambda} =\frac{1}{2}\|\cdot\|^2_2$, $A_1 = A$, and $A_2=B$, then \eqref{eq:primal damped second order flows} and \eqref{eq:multiplier damped second order flows} can recover the \eqref{eq: admm flows} by letting 
\begin{equation}
    \alpha_t = e^{-a_t}, \quad \beta_t = e^{2a_t + b_t}, \quad \gamma_t = e^{a_t} - \dot a_t + \zeta e^{a_t}. 
\end{equation}

We will provide a major extension to the \emph{constrained} accelerated distribution flows based on the proposed work and \cite{amir_BD_flow}.

\section{ACCELERATED CONSTRAINED DISTRIBUTION FLOWS}
To make a concrete connection between our proposed viewpoint and \cite{amir_BD_flow}, let's
first consider the following unconstrained optimization problem over the \emph{probability distribution} of finite second moment $\rho(x) = Law(X)\in \mathcal{P}_{ac,2}(\bbR^n)$ or $X \sim \rho, \; X \in \bbR^n$. 
\begin{equation}
\min_{\rho} F(\rho)
\end{equation} where functional $F(\rho): \mathcal{P}_{ac,2}(\bbR^n)\to \bbR$ is \emph{displacement} \emph{displacement}\footnote{One can think displacement convexity as geodesic convexity along $W_2$ geodesic curve $C(t) = (tT(X) + (1-t)X)_{\#}\rho_0$ with $C(0) = \rho_0,\; C(1) = \rho_1$ for some $\rho_0,\;\rho_1$, and $T_{\#}(\rho_0) = \rho_1$.} convex, which holds when $F(\rho) = D_{KL}(\rho|\rho_{\infty})$ with log-concave $\rho_{\infty}$. By letting $\phi = \frac{1}{2} \|\cdot\|^2$ and using the Wasserstein metric, \cite{amir_BD_flow} proposed the Hamiltonian flows for random variable $X\sim \rho$ and multiplier $Y$ (momentum in \cite{BD_flow1}, and also see B.4 in \cite{BD_flow1} with $P = Y$ ) as:
\begin{subequations} \label{eq: unconstraint dist hamiltonian flows}
\begin{align}  
    \frac{dX}{dt} & =\frac{\partial H}{\partial Y}(X,Y) = e^{a_t - g_t}Y, \\
    \frac{dY}{dt} &= -\frac{\partial H}{\partial X}(X,Y)=  -e^{a_t + b_t + g_t} \nabla_{\rho} F(\rho)(X),
\end{align}
\end{subequations}
with the Hamiltonian for the random variable defined as eq 4. in \cite{amir_BD_flow} or eq B.23 in \cite{BD_flow1}, $\nabla_{\rho} F(\rho): \bbR^n \to \bbR^n$ is the Wasserstein gradient defined in \cite{amir_BD_flow}, and $\nabla_{\rho}F(\rho)(X) = \nabla_{X}\frac{\delta F}{\delta \rho}(\rho)(X)$ with $\frac{\delta F}{\delta \rho}$ being the distribution functional derivative of $F$ at $\rho$.

Instead of constructing Hamiltonian flows for $X$ and multiplier (momentum) $Y$, we can equivalently construct the primal $X$ and mirror $Z_X$ flows with $\phi = \frac{1}{2}\|\cdot\|^2$ as follows:
\begin{subequations} \label{eq: unconstraint dist primal mirror flows}
\begin{align} 
    \frac{dX}{dt} & = e^{a_t}(\nabla \phi^* (U_X) - X) = e^{a_t}(U_X - X), \\
    \frac{dU}{dt} &=  -e^{a_t + b_t } \nabla_{\rho} F(\rho)(X),
\end{align}
\end{subequations} 
due to the self-duality of $\phi$, and $U_X = Z_X = X + e^{-a_t}\dot X$ is the lookahead term similar to the vector case. From \eqref{eq: unconstraint dist hamiltonian flows} and \eqref{eq: unconstraint dist primal mirror flows}, one can connect two formulations through $U_X = X + e^{-g_t}Y$. Thus, the Lagrangian and Hamiltonian flow formulations dual to each other, and this transformation is called the Legendre transform \cite{calculus_of_variatino_optimal_control}.

Define the probability distribution $\rho_i = Law(X_i)$ with $X_i\sim \rho_i,\; X_i \in \bbR^n$ being the random variable associated with $\rho_i$. 
We consider the constrained distributional optimization problem
\begin{equation} \label{problem: constrained distribution}
    \min_{\rho_1,\ldots,\rho_N} \Ssum F_i(\rho_i) \quad
    \text{s.t.} \quad \Ssum \mathcal{A}_i\rho_i = C,
\end{equation}
where $F_i(\rho_i): \mathcal{P}_{ac,2}(\bbR^n)\to \bbR$ is displacement convex and sufficiently smooth, $\mathcal{A}_i\rho_i = \int_{\bbR^n} \varphi_i(x)\rho_i(x)dx = \bbE_{\rho_i}[\varphi_i(X_i)]$ is linear in $\rho_i$ and $C \in \bbR^m$. 
Define the augmented Lagrangian and let $\phi_i(\cdot) = \Phi_{\Lambda}(\cdot) = \frac{1}{2}\|\cdot\|^2$ for \eqref{problem: constrained distribution} by 
\begin{equation}
    \mathcal{L}_{\mu}(\rho,\Lambda) = \Ssum F_i(\rho_i) + \langle\Lambda,r(\rho)\rangle + \frac{\mu}{2} \|r(\rho)\|^2
\end{equation} where $\rho = (\rho_1,\ldots,\rho_N)$, $\Lambda \in \bbR^m$ is the multiplier and $r(\rho) = \Ssum \mathcal{A}_i \rho_i- C$. We then assume \eqref{problem: constrained distribution} admits an optimal and saddle pair $(\rho^*,\Lambda^*)$. 
Consider the augmented BD Lagrangians:
\begin{subequations}
\begin{align}
    \bbL_{X}(X,\dot{X},t) &= e^{a_t + g_t}[\frac{1}{2}e^{-2a_t}\Ssum\bbE[\|\dot{X}_i\|^2 ] - e^{b_t}\mathcal{L}_{\mu}(\rho,U_{\Lambda})], \\
    \bbL_{\Lambda}(\Lambda,\dot{\Lambda},t) &= e^{a_t + g_t}[\frac{1}{2}e^{-2a_t}\|\dot{\Lambda}_i\|^2  + e^{b_t}\mathcal{L}_{\mu}(Z_{\rho},\Lambda)],
\end{align}
\end{subequations}
where $U_{\Lambda} = \Lambda + e^{-a_t}\dot \Lambda$, $Z_{\rho_i} = Law(U_i)$ and $Z_{\rho} = (Z_{\rho_1},\ldots,Z_{\rho_N})$ is the lookahead probability distribution with $U_i = Z_i = X_i + e^{-a_t}\dot{X}_i$ (due to the self duality of $\phi_i,\;\Phi_{\Lambda}$).

We can construct ODE flows similar to the vector case in \S\ref{sec: primal dual flows} as:
\begin{subequations} \label{eq: dist primal undamped flows}
\begin{align}  \label{eq: dist primal undamped flows: position}
&\dot X_i = e^{a_t}[U_i - X_i],\\ 
\label{eq: dist primal undamped flows: mirror}
&\dot U_i = -e^{a_t+b_t}  \nabla_{\rho_i}\mathcal{L_{\mu}}(\rho, U_{\Lambda}) (X_i),\\
\label{eq: dist multiplier undamped flow: position}
&\dot \Lambda = e^{a_t}[U_{\Lambda} - \Lambda],\\
\label{eq: dist multiplier undamped flow: mirror}
&\dot U_{\Lambda} = e^{a_t+b_t}  \nabla_{\Lambda}\mathcal{L_{\mu}}(Z_{\rho}, \Lambda)- \zeta \nabla^2 h_{\Lambda}(\Lambda)\dot \Lambda,
\end{align}
\end{subequations}  where $\nabla_{\rho_i}\mathcal{L_{\mu}}(\rho, \Lambda): \bbR^n \to \bbR^n$ is the Wasserstein gradient \cite{amir_BD_flow} for augmented Lagrangian. By the definition of the Wasserstein gradient, 
\[
\nabla_{\rho_i}\mathcal{L}_{\mu}(\rho,\Lambda) := \nabla \frac{\delta}{\delta \rho_i}\mathcal{L}_{\mu}(\rho,\Lambda),
\] with $\frac{\delta}{\delta \rho_i}\mathcal{L}_{\mu}(\rho,\Lambda)$ being the functional derivative of Lagrangian at $(\rho, \Lambda)$. See Appendix C of \cite{amir_BD_flow} for more information.
The next theorem states that the resulting first-order ODE flows can solve \eqref{problem: constrained distribution} with the exponential dual gap closing rate.

\begin{theorem}[Distributional convergence rates] \label{thm: exp duality gap dist}
Define $F(\rho) = \Ssum F_i(\rho_i)$ and assume \eqref{problem: constrained distribution} admitting to an optimal feasible saddle pair $(\rho^*,\Lambda^*)$. Under the transport-map condition  $\bbE \langle U_i  -T_{i,t}(X_i), \frac{d}{dt}T_{i,t}(X_i)\rangle = 0,\forall i \in \mathcal{N}$ in \cite{amir_BD_flow} with optimal quadratic push forward $(T_{i,t})_{\#}\rho_i = \rho_i^*$,
let $\mathcal{A}_i\rho_i = A_i\bbE[X_i]$, $A_i \in \bbR^{m\times n}$, the duality gap
\begin{equation} \label{eq: dist duality gap}
    G_{\mu}(\rho, \Lambda) = \mathcal{L}_{\mu}(\rho,\Lambda^*) - \mathcal{L}_{\mu}(\rho^*,\Lambda) = \mathcal{L}_{\mu}(\rho,\Lambda^*) - F(\rho^*),
\end{equation}
satisfies $G_{\mu}(\rho,\Lambda)=\mathcal{O}(e^{-b_t})$. Consequently, $\|r(\rho)\|^2,\,|F(\rho)-F(\rho^*)|^2=\mathcal{O}(e^{-b_t})$. 
Furthermore, if $\|U_{\Lambda}+\zeta\nabla h_{\Lambda}(\Lambda)\|\leq K<\infty$, then we have ``full'' rates $\|r(\rho)\|,\,|F(\rho)-F(\rho^*)|=\mathcal{O}(e^{-b_t})$, and a sharper energy dissipation bound is obtained if $\zeta>0$.

\begin{proof}
The proof is given in \ref{proof to dist convergence thm}.
\end{proof}
    
\end{theorem}
The constraint affects only the optimal distribution $\rho^*$ and doesn't alter the assumed properties of transport mappings; the condition holds automatically when $n=1$ \cite{amir_BD_flow}.

\section{EXAMPLES AND SIMULATIONS}
Two numerical examples, one considering fully distributed optimization over connected graphs and the other the constrained two-distribution transport problem, are provided to demonstrate the wide applicability of such methods.

\subsection{Fully distributed optimization}
\label{subsec: example 1}
Let $\mathcal{G(\mathcal{V},\mathcal{E})}$ represents a connected communication graph over vertices $\mathcal{V}$ and edges $\mathcal{E}$, the following is the standard distributed optimization problem in \cite{admm}:

\begin{equation}
\begin{aligned}
        \min_{(x_i,z_i\in \Delta_n),\forall i\in\mathcal{V}} &\Ssum \log(1+ \exp(-a_i^\top x_i))\\ 
        \text{s.t. } &x_i = L_iz_i,\quad\forall i \in \mathcal{V}|\lambda_i,\\ 
        &z_i = z_j, \quad\forall (i,j) \in \mathcal{E}|\theta_i ,
\end{aligned}
\end{equation} where $\Delta_n:=\{p\in\bbR^n_+|\|p\|_1 = 1\}$ is the probability simplex, and we choose $n=4$ with $a_i = [i,i/2,i+1,i+2]^\top$, $z_i\in \bbR^{P}$ is the local estimation of the global variable, $L_i \in \bbR^{n \times P}$ is the information selecting channel with $L_i = I$ and $P=4$ in this case, and $\lambda_i,\; \theta_i$ are the local multipliers. We also denote $\bar{z} = \frac{\Ssum z_i}{N}$ and $\eta_i = z_i - \bar{z}$ is the disagreement between agents (consensus is reached when $\|\eta_i\| = 0,\;\forall i$). 
The resulting first-order positional ODEs take the form of:
\begin{subequations}
\begin{align}
\dot x_i&=e^{a_t}\bigl(\nabla\phi_i^*(u_i)-x_i\bigr), \\
\dot z_i&=e^{a_t}\bigl(\nabla\phi_i^*(v_i)-z_i\bigr),\\
\dot\lambda_i
&=e^{a_t}\bigl(\nabla\varphi_i^*(w_i^\lambda)-\lambda_i\bigr),\\
\dot\theta_i
&=e^{a_t}\bigl(\nabla\varphi_i^*(w_i^\theta)-\theta_i\bigr),
\end{align}
\end{subequations} with strictly convex distance generating functions $\phi_i$ and $\varphi_i$ for the primal variables and the multipliers, respectively. In our simulation case, $\phi(x_i)  = \sum_{j=1}^n x_{i,j} \log x_{i,j} + I_{\Delta_n}(x_i)$ is the negative entropy generator plus indicator function of $\Delta_N$, and $\varphi_i = \frac{1}{2}\|\cdot\|^2$, and we choose $h_i(\cdot) = h_{\lambda}(\cdot) = \frac{1}{2}\|\cdot\|^2$. 

\begin{figure}
\centering
\includegraphics[width=0.4\textwidth]{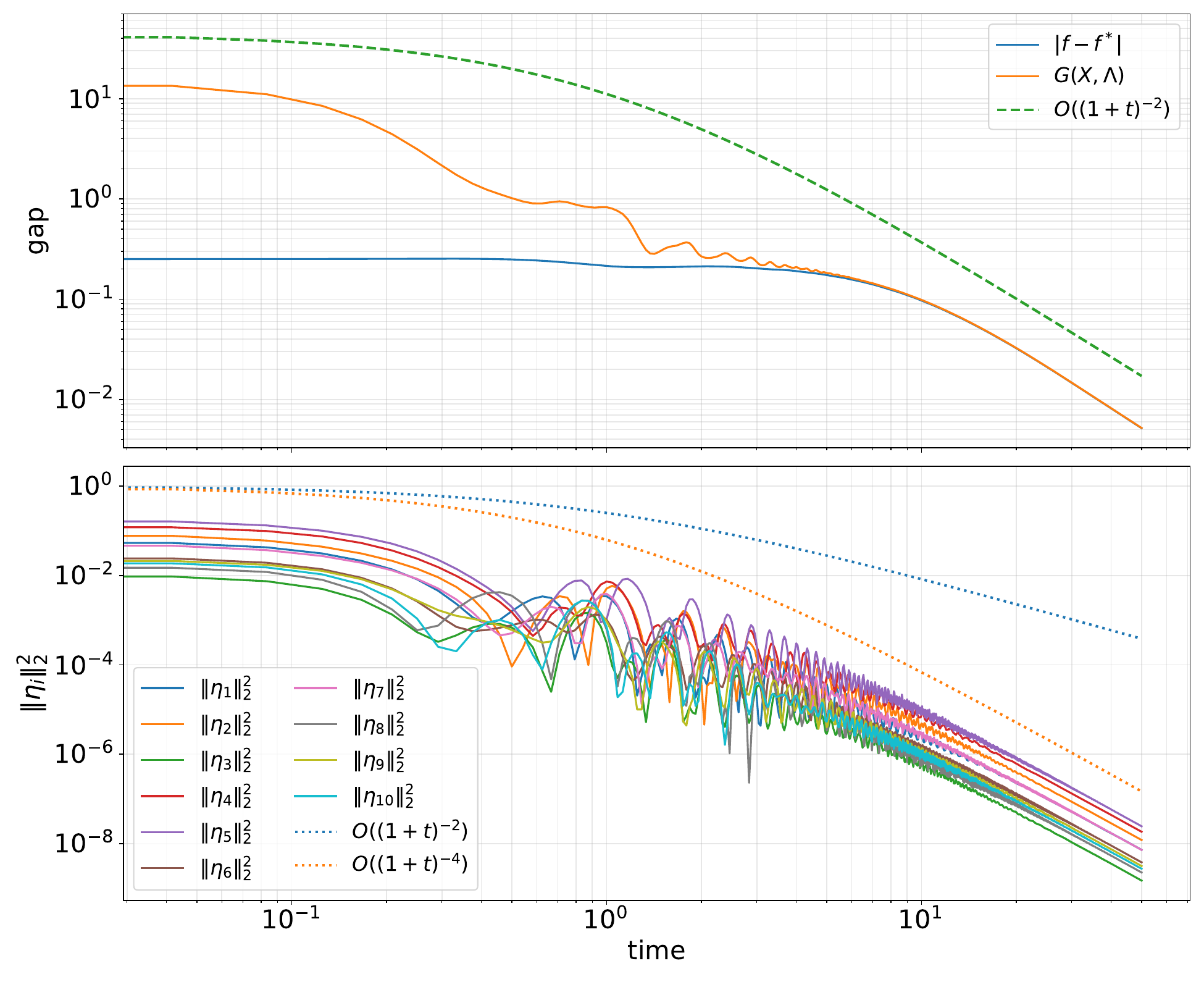}
\caption{Objective and duality gaps (upper) and consensus errors (lower)}
\label{fig: example 1}
\end{figure}

The corresponding mirror flows according to the general formulation \eqref{eq: damped mirror flows} are:
\begin{subequations}
\begin{align}
\dot u_i&=-e^{a_t+b_t}
\Bigl[\nabla f_i(x_i)+z_{\lambda_i}+\mu(x_i-L_i z_i)\Bigr]
-\zeta\dot x_i, \\
\dot v_i
&=- e^{a_t+b_t} \{- L_i^\top
\Bigl[z_{\lambda_i}+\mu(x_i-L_i z_i)\Bigr]
\nonumber\\
&\sum_{j\in\mathcal{N}_i}
\Bigl[(z_{\theta_i}-z_{\theta_j})+\mu(z_i-z_j)\Bigr]\}
-\zeta\dot z_i, \\
\dot w_i^\lambda
&=e^{a_t+b_t}(z_{x_i}-L_iz_{z_i})-\zeta\dot\lambda_i,\\
\dot w_i^\theta
&=e^{a_t+b_t}\sum_{j\in\mathcal{N}_i}(z_{z_i}-z_{z_j})-\zeta\dot\theta_i,
\end{align}
\end{subequations} with lookahead quantities $z_q = q + e^{-a_t}\dot q$ for $q = x_i, z_i, \lambda_i,\theta_i$.
The numerical simulation uses the complete graph with ten agents and adaptive numerical integrators, and the time scaling functions are chosen as 
\begin{equation}
    a_t = \log 2 - \log (1+t),\;b_t = 2\log(1+t),
\end{equation} with $\mu = 4$ and $\zeta = 1$. 
In the Fig. \ref{fig: example 1}, we denote $X = (x_1,z_1,\ldots,x_N,z_N)$, $\Lambda = (\lambda_1, \theta_i,\ldots,\lambda_N,\theta_N)$, and $f^*$ is the optimal value.

\subsection{NUMERICAL EXAMPLE: Constrained distributional optimization}
\begin{figure}
\centering
\includegraphics[width=0.49\textwidth]{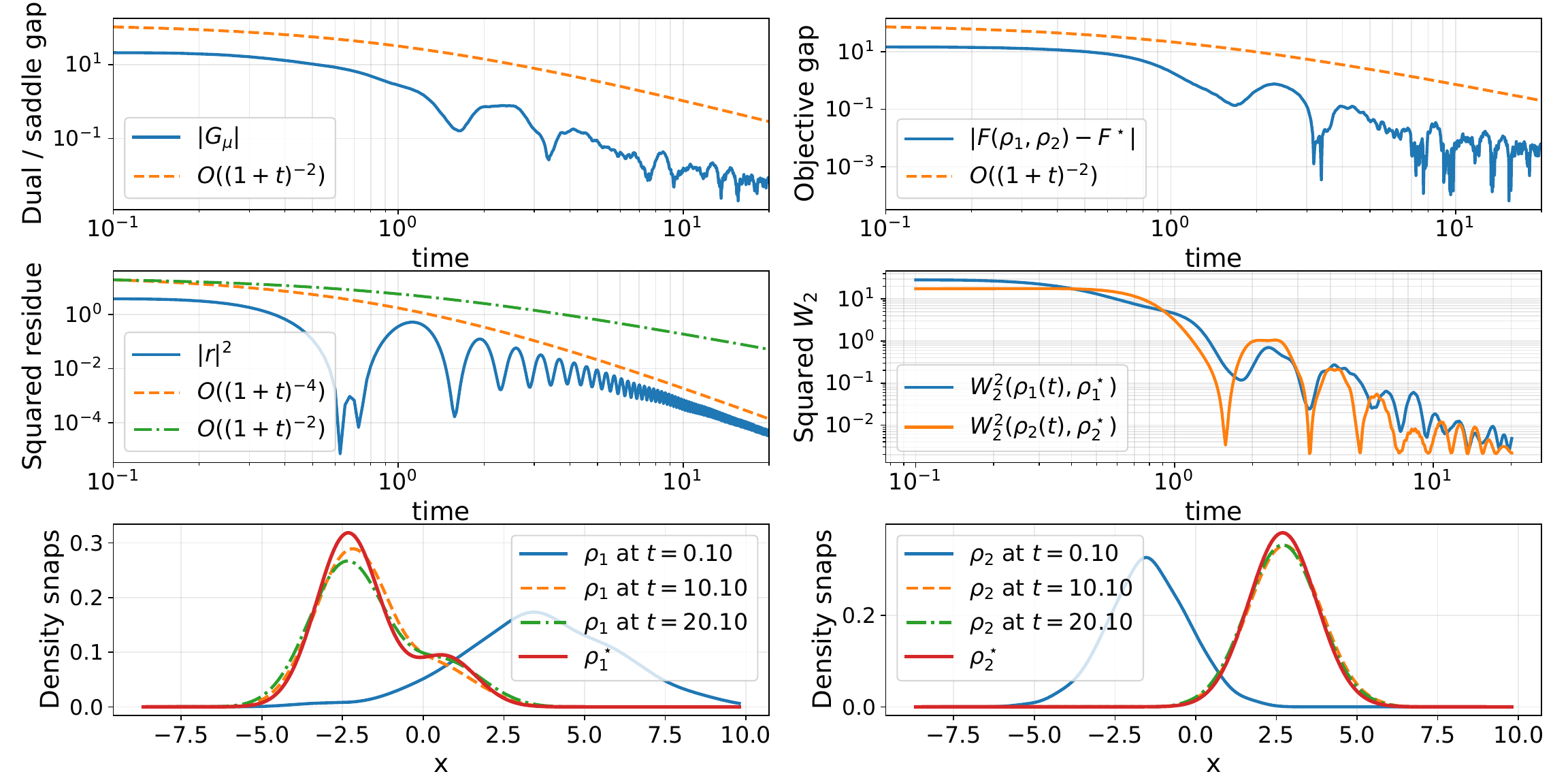}
\includegraphics[width=0.49\textwidth]{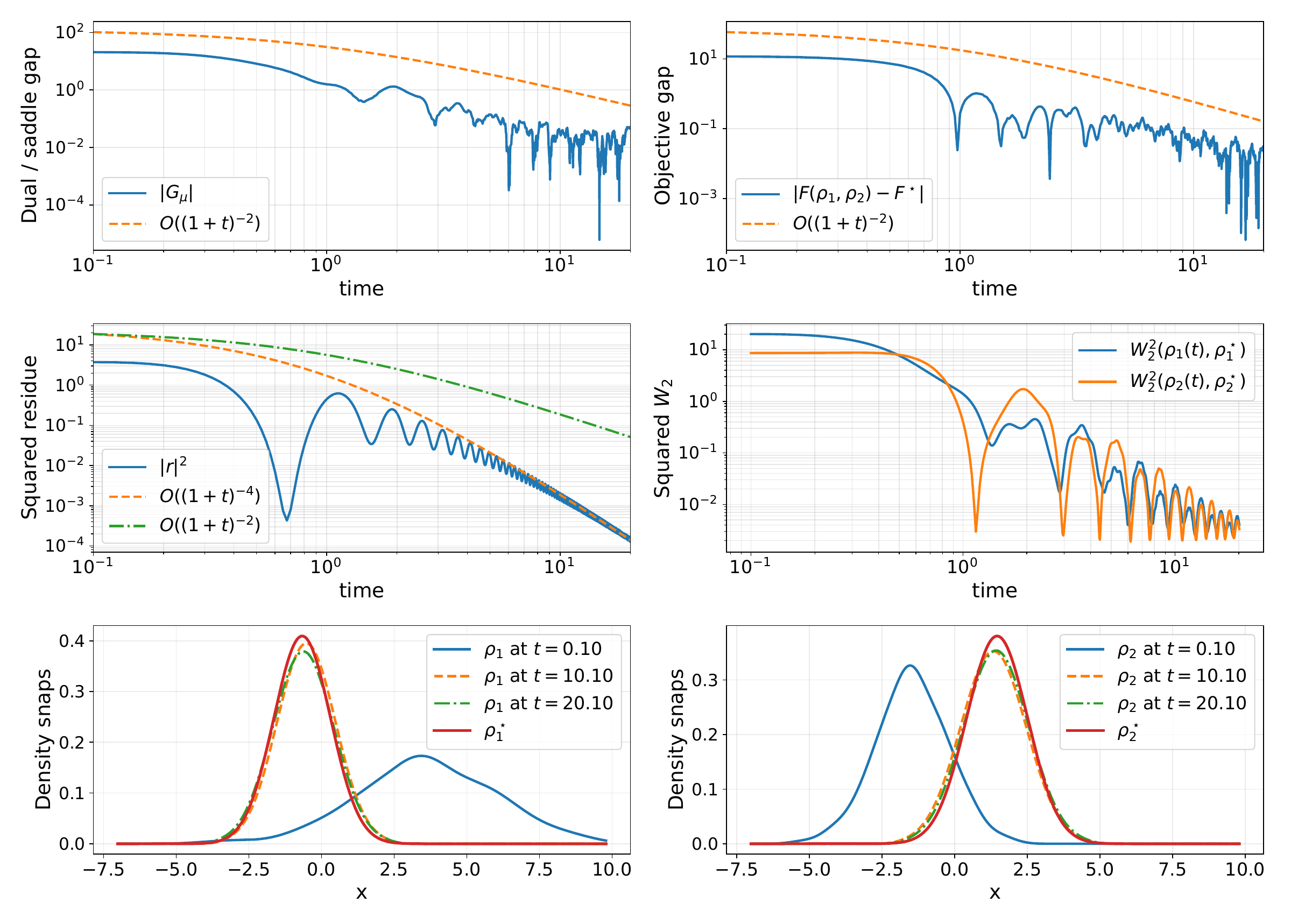}
\caption{Convergence properties and snapshots of the distributions for the bimodal target distribution (upper six) and normal distributions (lower six).}
\label{fig: distributional optimization}
\end{figure}

Consider the problem of minimizing the discrepancy between the distribution $\rho_i$ and the target distribution $\rho_{\infty}^i$ subject to an affine distributional constraint:
\begin{equation} 
    \min_{\rho_1,\rho_2}\sum _{i=1}^2 D_{KL}(\rho_i,\rho_{\infty}^i) \quad
    \text{s.t.} \quad \sum_{i=1}^2 \mathcal{A}_i\rho_i = C,
\end{equation} where $\rho_{\infty} = (\rho_{\infty}^1, \rho_{\infty}^2)$ is the target distribution with $\rho_{\infty}^1 = 0.55\mathcal{N}(-2,0.95) + 0.45\mathcal{N}(1,0.95)$ for bimodal case, $\rho_{\infty}^1 = \mathcal{N}(1,0.95)$ for the normal case, and $\rho_{\infty}^2 = \mathcal{N}(3,1.1)$,  $\sum _{i=1}^{2}\mathcal{A}_i\rho_i = \sum _{i=1}^{2}A_i\bbE[X_i] = C$ is the linear 1D distributional constraint with $A_1 = 1,\; A_2 = 0.8$, and $C = 0.5$.
For the KL divergence, the Wasserstein gradient takes the form of $\nabla_{\rho_i}D_{KL}(\rho_i|\rho_{\infty}^i)(X_i) = \nabla\log(\rho_i)(X_i) - \nabla\log (\rho_{\infty}^i)(X_i)$, and $\nabla \log(\rho_i)$ can be approximated by the finite amount of particles in Eqs. (23) and (24) in \cite{amir_BD_flow} for the normal and general distributions, respectively. Similarly, the distributions are approximated using finite particles, and the approximated flows yield the following ODEs when $h_{\Lambda} = \frac{1}{2}\|\Lambda\|^2$:
\begin{equation} 
\begin{aligned}
&\dot X_i^k = e^{a_t}[ U_i^k - X_i^k],\\ 
&\dot U_i^k = -e^{a_t+b_t}  [\tilde{ \nabla}_{\rho_i}D_{KL}(\rho_i |\rho_{\infty}^i)(X_i^k) + A_i^\top(U_{\Lambda} + \mu r(\rho))] \\
&\dot \Lambda = e^{a_t}[U_{\Lambda} - \Lambda],\dot U_{\Lambda} = e^{a_t+b_t} r(Z_{\rho}) - \zeta \dot \Lambda,
\end{aligned}
\end{equation} with particle index $k = 1,\ldots,200$, and $\tilde{ \nabla}_{\rho_i}D_{KL}(\rho_i |\rho_{\infty}^i)(X_i^k)$ the approximated Wasserstein gradient by using the diffusion map approximation \cite{amir_BD_flow}. We set, $\mu =3$, $\zeta = 0.8$, $a_t = \log 2 - \log (1+t),$ and $\;b_t = 2\log(1+t)$. Fig. \ref{fig: distributional optimization} empirically shows the $\mathcal{O}(e^{-b_t})$ rate and distribution transportation process over time, \emph{even though $D_{KL}(\rho_1|\rho^1_{\infty})$ is not displacement convex for the bimodal case.} The code can be found in \href{https://github.com/Justin900308/LCSS_opt_code}{Demo Code}.

\section{CONCLUSIONS}

By constructing augmented Bregman Lagrangians within an Euler–Lagrange–Rayleigh framework, we establish a unified variational construction for a class of accelerated primal–dual flows and recover representative existing methods as special cases. This viewpoint enables general Bregman geometries and, importantly, extends accelerated dynamics to affinely constrained optimization over probability distributions with $\mathcal{O}(e^{-b_t})$ convergence guarantees. More general convex inequality sets and non-smooth settings require future work involving the subdifferential and variational inequality constructions.


\section*{APPENDIX}
\subsection{Proof to Lemma \ref{lemma:Undamped ODES}} \label{proof: lemma Undamped}
Consider the variable $q$ that either represents $x_i$ or $\lambda$, with the following augmented BD Lagrangian for $q$:
\begin{equation}
    \bbL_{q}(q,\dot q,t) = e^{a_t+g_t}[\underbrace{D_{\phi_q}(z_q,q)}_{\text{Kinetic term}} -\text{sign} \underbrace{e^{b_t}V(q)}_{\text{Potential term}}], 
\end{equation} where $\text{sign} = 1$ if $q = x_i$ and $\text{sign} = -1$ if $q = \lambda$, and $V(q)$ is the variable dependent potential:
\begin{equation}
V(q) = 
\begin{cases}
    \mathcal{L}_{\mu}(\bx,z_{\lambda}),\quad q = x_i,\\
    \mathcal{L}_{\mu}(z_{\bx},\lambda),\quad q = \lambda.\
\end{cases}
\end{equation}
By \eqref{eq: euler lagrangian}, taking the derivatives for both terms leads to 
\begin{equation}
\begin{aligned} \label{eq: velocity derivative}
    &\frac{d}{dt}\frac{\partial \bbL_q}{\partial \dot q} \\&= e^{g_t}[\dot g(\nabla \phi_q(z_q) - \nabla \phi_q (q)) + \nabla^2\phi_q(z_q)\dot{z}_q - \nabla^2\phi(q)\dot q],
\end{aligned}   
\end{equation} and
\begin{equation}
\begin{aligned} \label{eq: positional derivative}
    \frac{\partial\bbL_q}{\partial q} &= e^{a_t + g_t}(\nabla \phi_q (z_q)-\nabla \phi_q(q)) - e^{g_t}\nabla^2 \phi_q(q)\dot q\\
    &-\text{sign }e^{a_t + b_t + g_t} \nabla_q V(q). 
\end{aligned}
\end{equation}
Subtracting \eqref{eq: velocity derivative} from \eqref{eq: positional derivative} and use the ideal scaling \cite{BD_flow1}, we arrive at 
\begin{equation} \label{eq: subtraction}
    \nabla^2\phi_q(z_q) \dot{z}_q + \text{sign } e^{a_t+b_t} \nabla_qV(q) = 0. 
\end{equation} Take the time derivative of $z_q$ and substitute into \eqref{eq: subtraction} with the assumption that $\nabla^2 \phi_q(z_q)$ is invertible, the final second-order ODE produces
\begin{equation}
    \ddot q + (e^{a_t} - \dot a_t)\dot q + \text{sign } e^{2a_t + b_t}[\nabla^2 \phi_q(z_q)]^{-1}\nabla_qV(q) = 0,  
\end{equation} and this recovered \eqref{eq: undamped second order flows} exactly for $q = x_i$ and $q= \lambda$. 

To show the equivalent corresponding first-order ODEs in \ref{subsec: undamped flows}, we aim to show that the positional and mirror flows together form the second-order ODEs. Recall the definition of the mirror flows with mirror coordinate $u_q = \nabla \phi_q(z_q)$:
\begin{equation}
    \dot u_q = - \text{sign } e^{a_t + b_t} \nabla_q V(q), 
\end{equation}
equal this with the direct time derivative of $u_q$
\begin{equation}
    \dot u_q = e^{-a_t} \nabla^2 \phi_q(z_q)[\ddot q + (e^{a_t} - \dot a_t)\dot q],
\end{equation} then the second-order ODE can be fully recovered from the first-order positional and mirror ODEs.

\subsection{Proof to Theorem \ref{theorem: convergence}}
\label{proof: convergence theorem}
Let $H_q = \nabla^2 h_q(q)$ and 
recall the Lyapunov function \eqref{eq:energy}, taking the time derivatives of each term gives:
\begin{equation}
\begin{aligned}
    \frac{d\mathcal{E}}{dt} = &\underbrace{\rate\Ssum[\langle \nabla_{x_i}\mathcal{L}_{\mu}(\bx,z_{\lambda}),x^*_i - z_i \rangle+\zeta\langle H_i\dot x_i,x_i^* - z_i \rangle]}_{\text{1.}}  \\
    +&\underbrace{\rate \langle \nabla_{\lambda}\mathcal{L}_{\mu}(z_{\bx},\lambda),z_{\lambda} - \lambda^*\rangle + \zeta\langle H_{\lambda}\dot \lambda, \lambda^* - z_{\lambda}\rangle}_{2.}\\
     + &\underbrace{[\rate \Ssum \langle \nabla_{x_i}\mathcal{L}_{\mu}(\bx,\lambda^*),z_i - x_i \rangle + \dot b_te^{b_t}G_{\mu}(\bx,\lambda))}_{3.} \\
     + & \underbrace{ [\Ssum \zeta \langle x_i - x_i^* ,H_i\dot x_i\rangle] 
     + \zeta \langle \lambda - \lambda^* ,H_{\lambda}\dot \lambda\rangle.}_{4.} 
\end{aligned}
\end{equation}
Group the damping terms together; one obtains 
\begin{equation} \label{eq: damping terms}
    \zeta [\langle H_q\dot q,q^* - z_q\rangle+\langle q-q^*,H_q \dot q\rangle ] = -\zeta e^{-a_t}\|\dot q\|^2_{H_q},
\end{equation}
since $z_q - q = e^{-a_t}\dot q$, and $q$ represents either $x_i$ or $\lambda$.  The next lemma relates the inner products to the dual gap:
\begin{lemma} [Inner products cancellation and inequality] \label{lemma: inner product cancellation}
\begin{equation}
\begin{aligned} \label{eq: coupling cancellation}
    &\Ssum [\langle \nabla_{x_i}\mathcal{L}_{\mu}(\bx,z_{\lambda}),x^*_i - z_i\rangle+
    \langle \nabla_{x_i}\mathcal{L}_{\mu}(\bx,\lambda^*),z_i - x_i\rangle]
    \\&+\langle \nabla_{\lambda}\mathcal{L}_{\mu}(z_{\bx},\lambda),z_{\lambda} - \lambda^* \rangle =  
    \Ssum \langle \nabla_{x_i}\mathcal{L}_{\mu}(\bx,\lambda^*),x^*_i - x_i\rangle \\
    &\leq
    \mathcal{L}_{\mu}(\bx^*,\lambda^*)
    -\mathcal{L}_\mu(\bx,\lambda^*)
    =-G_{\mu}(\bx,\lambda),
\end{aligned}
\end{equation}
\begin{proof}
    Let $\eta = z_{\lambda} - \lambda^*$
    By the construction of $\mathcal{L}_{\mu}$, $\nabla_{x_i}\mathcal{L}_{\mu}(\bx, z_{\lambda}) - \nabla_{x_i}\mathcal{L}_{\mu}(\bx, \lambda^*) = A_i^\top\eta $. Therefore, the difference between the left and right sides of the first equality of \eqref{eq: coupling cancellation} becomes 
    \begin{equation}
        \Ssum\langle A^\top_i \eta,x^*_i - z_i \rangle + \langle \eta, r(z_{\bx}) \rangle = \langle\eta, \Ssum A_i x^*_i - C \rangle = 0,
    \end{equation}
    and the last inequality of \eqref{eq: coupling cancellation} follows from convexity. 
\end{proof}

\end{lemma}
Note that the choice of $\phi$ (BD geometry) doesn't affect the above cancellation, so Lemma \ref{lemma: inner product cancellation} applies to different BD geometries. The lookahead term is also crucial for ensuring \eqref{eq:energy} produces the desired coupling rate \eqref{eq: coupling cancellation}.

Therefore, from \eqref{eq: damping terms}, \eqref{eq: coupling cancellation}, ideal scaling in \cite{BD_flow1}, and the positivity of $G_{\mu}(\bx,\lambda)$, we conclude
\begin{equation}
\begin{aligned}
    \dot{\mathcal{E}}\leq e^{b_t}(\dot b_t - e^{a_t})G_{\mu}(\bx,\lambda) - \zeta e^{-a_t}[\Ssum \|\dot x_i\|^2 + \|\dot \lambda\|^2] \leq 0.
\end{aligned}
\end{equation}
Lastly, by the construction of the Lyapunov function, $\mathcal{E} \geq e^{b_t}G_{\mu}(\bx,\lambda)$ together with  $G_{\mu}(\bx,\lambda) \geq \frac{\mu}{2}\|r(\bx)\|^2$, we have 
\begin{subequations}
\begin{align}
    \label{eq: gap convergence}
    e^{b_t}G_{\mu}(\bx,\lambda) &\leq \mathcal{E} \leq \mathcal{E}_0, \\
    \label{eq: residue reduction}
    \|r(\bx)\|^2 &\leq \frac{2\mathcal{E}_0}{\mu}e^{-{b_t}}. 
\end{align}
\end{subequations}
with initial total energy $\mathcal{E}_0$. 
Thus, we have the dual gap and residue rates: $G_{\mu}(\bx,\lambda) = \mathcal{O}(e^{-b_t})$ and $\|r(\bx)\|^2 = \mathcal{O}(e^{-b_t})$.
Lastly, $|F(\bx) - F(\bx^*)|^2 \leq [G_{\mu}(\bx,\lambda) + \|\lambda^*\|\|r(\bx)\|]^2$ implies 
\begin{equation} \label{eq: objective reduction}
    |F(\bx) - F(\bx^*)|^2 = \mathcal{O}(e^{-b_t}), 
\end{equation} by the gap and residue reduction rates in \eqref{eq: gap convergence} and \eqref{eq: residue reduction}.

\begin{lemma}\label{lemma: boundedness}
Recall the strongly convexity of $\Phi_{\lambda}$ and $h_{\lambda}$, we state there exists some $K <\infty$ such that $\|u_{\lambda} + \zeta\nabla h_{\lambda}(\lambda)\|\leq K$.

\begin{proof}
    To bound $\|u_{\lambda} + \zeta\nabla_{\lambda}h(\lambda)\|$, it is sufficient to show $\|u_{\lambda}\|\leq K_{\lambda}<\infty$ and $\|\nabla h_{\lambda}(\lambda)\|\leq K_h < \infty$. 
    By the strong convexity of $\Phi_{\lambda}$, there is $m_{\Phi}>0$ so that $\mathcal{E}_0\geq D_{\Phi_{\lambda}}(\lambda^*,z_{\lambda}) \geq m_{\Phi}\|\lambda^* - z_{\lambda}\|^2$, and the first inequality comes from the construction of the Lyapunov function. Similarly, we have $\mathcal{E}_0\geq \zeta D_{h_{\lambda}}(\lambda^*,\lambda) \geq \zeta m_{h}\|\lambda^* - \lambda\|^2$, with some $m_h>0$. Now, recall $u_{\lambda} = \nabla{\Phi_{\lambda}}(z_{\lambda})$, and both $\nabla\Phi_{\lambda}, \; \nabla h_{\lambda}$ are both continuous. By continuity and the boundedness of $z_{\lambda}$ and $\lambda$, there exist $K_{\lambda},K_h < \infty$ such that $\|u_{\lambda}\| = \|\nabla\Phi_{\lambda}(z_{\lambda})\| \leq K_{\lambda}$ and $ \|\nabla h_{\lambda}(\lambda)\| \leq K_h$; thus there is some $K < \infty$ and $\|u_{\lambda} + \zeta\nabla h(\lambda)\| \leq K$. 
\end{proof}
    
\end{lemma}

Now, we show that when $\zeta > 0$ and  $\|u_{\lambda} + \zeta \nabla h_{\lambda}(\lambda) \| \leq K < \infty$, the proposed flows also exhibit exponential convergence for the residual and objective. 

By the definition of residue and lookahead term, $r(z_{\bx}) = r(\bx) + e^{-a_t}\dot r(\bx)$, and recall the multiplier mirror flow in \eqref{eq: damped mirror flows}, we write
\begin{equation} \label{eq: new multiplier flow}
    \dot u_{\lambda} + \zeta \nabla^2 h(\lambda) \dot \lambda = e^{a_t + b_t}r(\bx) + e^{b_t}\dot r(\bx).
\end{equation}
Subtracting \eqref{eq: new multiplier flow} from $\frac{d}{dt}(e^{b_t}r) = e^{b_t}\dot r(\bx) + \dot b_te^{b_t}r(\bx)$ gives 
\begin{equation} \label{eq: zero rate}
\dot q_t=(\dot b_t-e^{a_t})(q_t+u_{\lambda}+\zeta\nabla h_{\lambda}(\lambda)),
\end{equation}
where $q_t=e^{b_t}r(\bx)-u_{\lambda}-\zeta\nabla h_{\lambda}(\lambda)$. By ideal scaling \eqref{eq: ideal scaling} and $\|u_{\lambda}+\zeta\nabla h_{\lambda}(\lambda)\|\leq K$, which gives $\|q_t\|\leq \max{\{\|q_0\|,K}\}<\infty$ from stable forced linear solution, and hence $\|r(\bx)\|=\mathcal{O}(e^{-b_t})$.
With the exponential rate for the residue, we also have $|F(\bx) - F(\bx^*)| = \mathcal{O}(e^{-b_t})$; thus, under the boundedness assumption, full rates are achieved.

\subsection{Proof to Theorem \ref{thm: exp duality gap dist}} \label{proof to dist convergence thm}
Define the optimal push forward $(T_{i,t})_{\#}\rho_i = \rho^*_i$. Consider the Lyapunov function (similar to \eqref{eq:energy} for the vector case with $\phi = \frac{1}{2}\|\cdot\|^2$ and $x^*_i = T_{i,t}(X_i)$)
\begin{equation} \label{eq: energy distribution}
\begin{aligned}
    \mathcal{E}_{\rho} &= \underbrace{[\frac{1}{2}\Ssum\bbE[\|U_i - T_{i,t}(X_i)\|^2]]}_{1.}+ \underbrace{\frac{1}{2}\|U_{\Lambda} - \Lambda^*\|^2}_{2.} \\
    &+ \underbrace{e^{b_t}G_{\mu}(\rho,\Lambda)}_{3.}+\underbrace{\zeta D_{h_{\Lambda}}(\Lambda^*, \Lambda)}_{4.}.
\end{aligned}
\end{equation} 
Note that with our choice of $\mathcal{A}_i$, the Wasserstein gradient terms for the augmented Lagrangian become
\begin{equation} \label{eq: energy distribution}
    \nabla_{\rho_i} \mathcal{L}_{\mu}(\rho,\Lambda)(X_i) = \nabla_{\rho_i}F_i(\rho_i)(X_i) +  A_i^\top (\Lambda + \mu r(\rho)),
\end{equation}
and \eqref{eq: energy distribution} produces the time derivative
\begin{equation}
\begin{aligned}
\dot{\mathcal{E}}_{\rho} &= 
\underbrace{\rate \Ssum\bbE[\langle \nabla_{\rho_i}\mathcal{L}_{\mu}(\rho,\Lambda)(X_i) ,U_i - T_{i,t}(X_i) \rangle]}_{1.} \\ 
&+ \underbrace{\rate \langle U_{\Lambda}-\Lambda^*,r(z_{\rho}) \rangle}_{2.} \\
&+\underbrace{\dot{\beta}_te^{b_t} + \rate \Ssum \bbE[\langle \nabla_{\rho_i}\mathcal{L}_{\mu}(\rho,\Lambda^*)(X_i), U_i - X_i \rangle]}_{3.}\\
&+ \underbrace{\zeta \langle H_{\Lambda}\dot{\Lambda},\Lambda^* - U_{\Lambda}\rangle}_{4.}.
\end{aligned}
\end{equation}

Thus, assume the technical assumption $\bbE \langle U_i  -T_{i,t}(X_i), \frac{d}{dt}T_{i,t}(X_i)\rangle = 0,\forall i \in \mathcal{N}$ \cite{amir_BD_flow}, collect the inner product terms and use Lemma \ref{lemma: inner product cancellation} with $\eta = U_{\Lambda} - \Lambda^*$, one can easily show  
\begin{equation}
\begin{aligned}
    \Ssum\bbE &\langle A_i^\top\eta, T_{i,t}(X_i) - U_i \rangle + \langle\eta, r(Z_{\rho})  \rangle \\
    & = \langle \eta, \Ssum \mathcal{A}_i \rho_i^* - C \rangle = 0,
\end{aligned}
\end{equation} since the $(T_{i,t})_{\#}\rho_i = \rho_i^*$ and optimal distribution $\rho^*$ is feasible.
By the displacement convexity in $\rho$, we have a similar inequality condition as in the vector case \ref{lemma: inner product cancellation}:
\begin{equation}
    \Ssum \bbE\langle \nabla_{\rho_i}\mathcal{L_{\mu}}(\rho,\Lambda^*)(X_i),T_{i,t}(X_i) - X_i\rangle \leq -G_{\mu}(\rho,\Lambda);
\end{equation} therefore:
\begin{equation}
    \dot{\mathcal{E}}_{\rho} \leq e^{b_t}(\dot b_t - e^{a_t})G_{\mu}(\rho,\Lambda) \leq 0,
\end{equation} and this implies $G_{\mu}(\rho, \Lambda) = \mathcal{O}(e^{-b_t})$ and $\|r(\rho)\|^2 = \mathcal{O}(e^{-b_t})$, since $G_{\mu}(\rho,\Lambda) \geq \frac{\mu}{2}\|r(\rho)\|^2$. Since the constraint and the multiplier are of finite dimension, the proof of the exponential residue decrease is largely identical to the proof made in \ref{proof: convergence theorem} and \ref{lemma: boundedness} for the vector case; hence, it is omitted. With $\|r(\rho)\| = \mathcal{O}(e^{-b_t})$, it is clear that $F(\rho) - F(\rho^*) = \mathcal{O}(e^{-b_t})$.


\bibliographystyle{ieeetr}
\bibliography{references}

\end{document}